\documentclass[11pt,final]{aims}
\usepackage{amsmath}
\usepackage[T1]{fontenc}
\usepackage[latin1]{inputenc}
\usepackage{times}
\usepackage{amsfonts}
\usepackage{fullpage}
\usepackage{amsmath}
  \usepackage{paralist}
\usepackage{treetex}

\usepackage{latexsym,amsfonts,amssymb,amsmath}
\usepackage{times}
\usepackage{paralist}
\usepackage{treetex}

\def\theequation{\thesection.\arabic{equation}}
\newcommand{\refs}[1]{{\rm(\ref{#1})}}

\newtheorem{theorem}{Theorem}[section]
\newtheorem{lemma}[theorem]{Lemma}

\newtheorem{remark}[theorem]{Remark}

\def\d{\hbox{d}}

\title{Asymptotic analysis of stochastic variational inequalities modeling an elasto-plastic problem with vanishing jumps\thanks{This research was partially supported by a grant from CEA, Commissariat \`a l'\'energie atomique and by the National Science Foundation under grant DMS-0705247}}
\date{\today}
\begin{document}\maketitle
\noindent
\author{{\sc Alain Bensoussan}\footnote{This research in the paper was supported by WCU (World Class University) program through the National Research Foundation of Korea funded by the Ministry of Education, Science and Technology (R31 - 20007).}\\[2pt]
International Center for Decision and Risk Analysis, School of Management, University of Texas at Dallas, Box 830688, Richardson, Texas 75083-0688,\\[3pt] 
Graduate School of Business, the Hong Kong Polytechnic University,Graduate Department of Financial Engineering, Ajou University.\\[3pt]
alain.bensoussan@utdallas.edu\\[6pt]
{\sc Hector Jasso-Fuentes}\\[2pt]
Departamento de Matem\'aticas, CINVESTAV-IPN, A, Postal 14-740, M\'exico D.F. 07000,\\[3pt]
M\'exico, hjasso@math.cinvestav.mx\\[6pt]
{\sc St\'ephane Menozzi}\\[2pt]
Universit\'e Denis Diderot-Paris 7, LPMA, 175 rue du Chevaleret 75013 Paris,\\[3pt]
menozzi@math.jussieu.fr\\[6pt]
{\sc Laurent Mertz}\\[2pt]
Laboratoire Jacques-Louis Lions,Universit\'e Pierre et Marie Curie, 4 place jussieu, Paris, 75005 France,\\[3pt]
mertz@ann.jussieu.fr
\maketitle
\begin{abstract}
{In a previous work by the first author with J. Turi \cite{BenTuri1}, a stochastic variational inequality has been introduced to model an elasto-plastic oscillator with noise. 
A major advantage of the stochastic variational inequality is to overcome the need to describe the trajectory by phases (elastic or plastic). This is useful, since the sequence of phases cannot be characterized easily. In particular, when a change of regime occurs, there are numerous small elastic phases which may appear as an artefact of the Wiener process. However, it remains important to have informations on both the elastic and plastic phases. In order to reconcile these contradictory issues, we introduce an approximation of stochastic variational inequalities by imposing artificial small jumps between phases allowing a clear separation of the elastic and plastic regimes.
In this work, we prove that the approximate solution converges on any finite time interval, when the size of jumps tends to $0$.} 
\end{abstract}

\section{Introduction}
The elastic-perfectly-plastic (EPP) oscillator under standard white noise excitation is the simplest structural model exhibiting a hysteretic  behavior. Moreover, the model is representative of the behavior of mechanical structures which vibrate mainly on their first deformation mode.  
In the context of earthquake engineering, relevant applications to piping systems under random vibrations can be accessed this way \cite{Feau1,Feau2}.
The main difficulty to study these systems comes from a frequent occurrence of nonlinear phases (plastic phases) on small time intervals. A nonlinear phase corresponds to a permanent deformation, or in other words to a plastic deformation. A plastic deformation is produced when the stress of the structure exceeds an elastic limit.
Denoting by $x(t)$ the elasto-plastic displacement, we consider the problem
\begin{equation}
\label{chap4:hysteresis1}
\ddot{x} + c_0 \dot{x} + \mathbf{F}(x(s), 0 \leq s \leq t) = \dot{w},
\end{equation}
with initial conditions of displacement and velocity
\[
x(0)=x \quad, \quad \dot{x}(0) = y.
\]
Here $c_0>0$ is the viscous damping coefficient, $k>0$ the stiffness, $w$ is a Wiener process and $\mathbf{F}(\{x(s), 0 \leq s \leq t\})$ is a nonlinear functional which depends on the entire trajectory $\{ x(s), 0 \leq s \leq t\}$ up to time $t$.
Denote $y(t):= \dot{x}(t)$. Equation \eqref{chap4:hysteresis1} written as a stochastic differential equation (SDE) reads
\begin{equation}\label{chap4:hysteresis1bis}
	\d y(t)= -(c_0 y(t) + F(\{ x(s), 0 \leq s \leq t \}))\d t + \d w(t), \quad d x(t) = y(t) \d t.
\end{equation}
Beyond a given threshold $\vert F(\{ x(s), 0 \leq s \leq t \}) \vert=kY$ for the nonlinear restoring force, the material goes through plastic deformation (see e.g. \cite{KarSchar}). Introducing $\Delta(t)$, the total plastic yielding accumulated up to time $t$, we can define a new state variable $z(t)$ as $z(t):=x(t)-\Delta(t)$. It follows that in the plastic regime, $\dot z(t)=0$.
From now on, we choose to express the restoring force $F(\{ x(s), 0 \leq s \leq t \})$ in \eqref{chap4:hysteresis1bis} in terms of the new variable $z(t)$ as $F(\{ x(s), 0 \leq s \leq t \}):=kz(t)$ where $ \vert z(t) \vert  \leq Y$. In other words we consider a linear restoring force of the variable $z(t)$ (whose modulus equals Y during the plastic phases). This type of force characterizes the elasto-perfectly-plastic behavior.

In \cite{BenTuri1}, for the previous choice for $F$, it is shown that \eqref{chap4:hysteresis1bis} is equivalent to a stochastic variational inequality (SVI). In addition, existence and uniqueness of an invariant measure for the solution of SVI have also been proven. For a general framework dealing with this class of inequalities we refer the reader to \cite{BenLio} and to \cite{DuvLion} for specific deterministic applications to mechanics.
Although SVIs have been already studied in \cite{BenLio} to represent reflection-diffusion processes in convex sets, no connection with random vibration problems had been made so far.
From \cite{BenTuri1}, the solution $(y(t),z(t)) \in \mathbb{R}^2$ of \eqref{chap4:hysteresis1bis} satisfies
\begin{equation}\label{chap4:SVI}
\dot{y}(t)= -(c_0 y(t) + k z(t)) + \dot{w}(t), \quad
(\dot{z}(t)-y(t)) (\phi - z(t)) \geq 0,\quad 
\forall  \vert \phi \vert \leq Y, \quad
\vert z(t) \vert \leq Y.
\end{equation}
In terms of dynamics of the process $(y(t),z(t))$, a plastic deformation begins when $z(t)$ reaches and is absorbed by $Y$ (resp. $-Y$) with positive (resp. negative) slope, $y(t) >0$. (resp. $y(t)<0$) i.e. when $\textup{sign}(y(t)) z(t) = Y$. Then, the plastic behavior ends when the velocity changes sign. At that time, the elastic behavior is reactivated. However, around 0, the velocity which is subjected to white noise, changes sign an infinite number of times during any small time interval. Often, this leads to a return into plastic behavior in a short time duration. 
This phenomenon is called \textit{micro-elastic phasing} and has been studied in \cite{FeauMer} using the numerical method developed in \cite{BenMerPiTu} for the SVI \eqref{chap4:SVI}. It plays a crucial role on frequency and statistics of plastic deformations. 
Because of this phenomenon, frequency of occurence, statistics (time duration or absolute plastic deformation) and the sequence of entry in plastic phase (as well as the sequence of exit) are not well defined.
In this paper, we consider an EPP oscillator under standard white noise excitation subjected to jumps (presented below) to study phase transitions. It has the advantage of separating phases clearly, while being an approximation. We prove the convergence of the approximated process towards the solution of the stochastic variational inequality \eqref{chap4:SVI}. 
 
\subsection{Model definition and convergence results}\label{chap4:sec2}

In this subsection, we introduce a stochastic variational inequality whose dynamics is ``almost" similar to the one of  \eqref{chap4:SVI} except that the second component is subjected  to jumps of magnitude $\varepsilon>0$ at some random times corresponding to the various exits of the plastic phases. 

Precisely, we describe the evolution of the new process $(y^\epsilon(t),z^\epsilon(t))$ by the following procedure; we start by defining $\tau_0^\epsilon := 0$ and by $(y^\epsilon_0(t),z^\epsilon_0(t))$ the solution of \eqref{chap4:SVI}, with initial conditions:
\[
y^\epsilon_0(0)= y \quad  \mbox{and} \quad  z^\epsilon_0(0)= z,\ (y,z)\in {\mathbb R}\times (-Y,Y):=D.
\]
Then, we define
\[
 \tau_1^\epsilon:=\inf\{t > 0, \quad y_0^\epsilon(t)=0 \quad \mbox{and} \quad \vert z_0^\epsilon(t) \vert=Y\}.
\]
For $t\geq\tau_1^\epsilon$, let $(y^\epsilon_1(t),z^\epsilon_1(t))$ be the solution of \eqref{chap4:SVI}
with initial conditions:
\[
y^\epsilon_1(\tau_1^\epsilon)=0 \quad \mbox{and} \quad z^\epsilon_1(\tau_1^\epsilon)=\textup{sign}(z^\epsilon_0(\tau^\epsilon_1)) \left(Y-\epsilon\right),
\]
again, we define
\[
\tau_{2}^\epsilon:=\inf\{t > \tau^\epsilon_1, \quad y_1^\epsilon(t)=0 \quad  \mbox{and} \quad \vert z_{1}^\epsilon(t) \vert =Y\}.
\]
In a recurrent manner, knowing $\tau_n^\epsilon$, $y_n^\epsilon(t)$, and $z_n^\epsilon(t)$, we define
\[
\tau_{n+1}^\epsilon:=\inf\{t > \tau^\epsilon_n, \quad y_n^\epsilon(t)=0 \quad  \mbox{and} \quad \vert z_{n}^\epsilon(t) \vert =Y\},
\]
and
$(y^\epsilon_{n+1}(t),z^\epsilon_{n+1}(t))$ be the solution of \eqref{chap4:SVI} with initial conditions:
\[
y^\epsilon_{n+1}(\tau^\epsilon_{n+1})=0 \quad \mbox{and} \quad z^\epsilon_{n+1}(\tau^\epsilon_{n+1})= \textup{sign}(z^\epsilon_n(\tau^\epsilon_{n+1})) \left(Y-\epsilon\right).
\]
Now, we define the process $(y^\epsilon(t),z^\epsilon(t))$ on each interval of time $[\tau^\epsilon_n,\tau^\epsilon_{n+1})$ as follows:

\begin{equation}\label{chap4:SVIe}
\dot{y}^\epsilon(t)= -(c_0 y^\epsilon(t) + k z^\epsilon(t)) + \dot{w}(t), \quad (\dot{z}^\epsilon(t)-y^\epsilon(t)) (\phi - z^\epsilon(t)) \geq 0, \quad \forall \vert \phi \vert \leq Y, \quad  \vert z^\epsilon(t) \vert  \leq Y
\end{equation}
with the following jump-conditions:
\[
y^\epsilon(\tau^\epsilon_{n}-)=0,\qquad z^\epsilon(\tau^\epsilon_{n}-)= z^\epsilon_{n-1}(\tau^\epsilon_{n}),
\]
and
\[
y^\epsilon(\tau^\epsilon_{n}) = 0,\qquad z^\epsilon(\tau^\epsilon_{n}) = \textup{sign}(z^\epsilon_{n-1}(\tau^\epsilon_{n}))(Y-\epsilon).
\]
\begin{remark}
\label{chap4:rem-reg}
By construction, the process $(y^\epsilon(t),z^\epsilon(t))$ is c\`adl\`{a}g; hence it is regular. In particular,  for each \underline{fixed} time $T>0$, the number of jumps arise in $(0,T]$, is finite a.s.
\end{remark}
We will prove that the solution $(y^\epsilon(t),z^\epsilon(t))$ converges to $(y(t),z(t))$ on any finite time interval, when $\epsilon$ goes to $0$ in the sense described below.

\section{Main results}
Our main result is the following theorem.
\begin{theorem}\label{chap4:thm}
Fix $T>0$, and consider the processes $(y(t),z(t))$ and $(y^\epsilon(t),z^\epsilon(t))$ satisfying \eqref{chap4:SVI} and \eqref{chap4:SVIe} respectively. Suppose that $k > X_+(c_0) :=\frac{1}{2} \left ( -\frac{c_0}{3} + c_0\sqrt{\frac{1}{9} + 4 \frac{c_0}{6}} \right )$. Then, the following convergence property holds:
\[
\frac{1}{\epsilon}\mathbb{E}\left[ \sup_{0 \leq t \leq T}\left\{ \vert y(t)-y^\epsilon(t) \vert^2+k\left \vert z(t)-z^\epsilon(t)\right \vert^2\right\}\right]\to 0 \quad \mbox{as} \quad \epsilon\textcolor{black}{\to} 0.
\]
\end{theorem}
\begin{remark}
Observe that the above condition relating $k$ and $c_0$ is purely technical. It will appear clearly in the proof of Lemma \ref{chap4:lemma2} below.
\end{remark}
\subsection{Preliminary results}
For $(y,z) \in D:=\mathbb{R} \times (-Y,Y)$, we consider the ``elastic" process $(y_{yz}(t),z_{yz}(t))$:
\begin{align*}
z_{yz}(t)  & = e^{\frac{-c_0t}{2}} \lbrace z \cos{(\omega t)} + \frac{1}{\omega}(y+\frac{c_0}{2}z) \sin{(\omega t)}\rbrace +  \frac{1}{\omega}  \int_0^t e^{-\frac{c_0}{2}(t-s)} \sin{(\omega(t-s))} \d w(s),\\
y_{yz}(t) & =-\frac{c_0}{2}z_{yz}(t)+ e^{-\frac{c_0 t}{2}} \lbrace -\omega z \sin{(\omega t)} + (y+\frac{c_0}{2}z)\cos{(\omega t)}\rbrace +\int_0^t e^{-\frac{c_0}{2}(t-s)}\cos{(\omega (t-s))} \d w(s).
\end{align*}
where, assuming $4k >c_0^2$,
\[
\omega := \frac{\sqrt{4k - c_0^2}}{2}.
\]

\begin{remark}
The terminology ``elastic'' is justified from the observation that $(y_{yz}(t),z_{yz}(t))$ is actually the solution of
\[
\dot{y}(t)= -(c_0 y(t) + k z(t)) + \dot{w}(t), \quad \dot{z}(t)=y(t), \quad (y_{yz}(0),z_{yz}(0))=(y,z),
\]
that is the explicit solution of \eqref{chap4:SVI} when the threshold $Y=\infty$ (purely elastic case). Note that the condition $4k >c_0^2$ is needed so that $(y(t),z(t))$ have real valued solutions.
\end{remark}

Define 
\begin{equation}\label{chap4:theta-eps}
\theta(y,z):= \inf \{ t>0, \quad \vert z_{yz}(t) \vert = Y\},
\end{equation}
where $(y_{yz}(0),z_{yz}(0))=(y,z)$. For $t \in [0,T]$, we set $u(y,z,t):=\mathbb{P}[\theta(y,z)>T-t]$. 
This function is regular and satisfies the mixed Cauchy-Dirichlet parabolic PDE
\begin{equation}
\label{chap4:pde1}
-u_t + Au = 0, \mbox{ in } D; 
\quad u(y,Y,t)= 0, \quad  y>0; 
\quad u(y,-Y,t)=0, \quad  y<0;
\quad u(y,z,T)= 1,
\end{equation}
with 
\[
Au=-\dfrac{1}{2} u_{yy} + (c_0y+kz)  u_{y} - y u_z.
\]
For $t<T$, the function $u(y,z,t)$ is locally smooth.
On the other hand, in the particular case when $(y_{yz}(0),z_{yz}(0)) := (0,Y-\epsilon)$, we consider the probability density function $p^\epsilon$ of $(y_{0,Y-\epsilon}(t),z_{0,Y-\epsilon}(t))$. It is also known that $p^\epsilon$ satisfies Chapman-Kolmogorov's equation
\begin{equation}\label{chap4:bound-cond-pe}
p^\epsilon_t + A^*p^\epsilon = 0, 
\quad p^\epsilon(y,z,0)= \delta_{0,Y-\epsilon}(y,z),
\end{equation}
where $A^*$ represents the adjoint operator of $A$; that is
\[
A^*p^\epsilon=-\frac{1}{2} p^\epsilon_{yy} - ((c_0y+kz)p^\epsilon)_y + y p^\epsilon_z.
\]
Next, observe that the processes $z_{0,Y-\epsilon}(t)$ and $y_{0,Y-\epsilon}(t)$ are gaussian processes. The key point is to express the solution of \eqref{chap4:pde1} through its variational formulation with $p^\epsilon$ as test function (see proof of Lemma \ref{chap4:lemma1}).\\

The mean, variance and covariance of $z_{0,Y-\epsilon}(t)$ and $y_{0,Y-\epsilon}(t)$ write:
\begin{equation}
\label{chap4:meanvarze}
m^\epsilon(t):=(Y-\epsilon)e^{-\frac{c_0t}{2}}(\cos \omega t+\frac{c_0}{2 \omega }\sin  \omega t ), \quad
\sigma_z^2(t):=\frac{1}{\omega^2}\int_0^t e^{-c_0s} \sin^2 (\omega s) \d s,
\end{equation}
\begin{equation}
\label{chap4:meanvary}
q^\epsilon(t):=-(Y-\epsilon)\frac{k}{\omega} e^{-\frac{c_0t}{2}} \sin \omega t , \quad
\sigma_y^2(t):=\int_0^t e^{-c_0s} (\cos \omega s-\frac{c_0}{2\omega}\sin \omega s)^2 \d s,
\end{equation}
and
\begin{equation}
\label{chap4:covyz}
\sigma_{yz}(t):=\frac{1}{2 \omega^2}e^{-c_0t} \sin^2 \omega t.
\end{equation}
The density $p^\epsilon$ then explicitly writes
\begin{align}
p^\epsilon(y,z,t)  = & \frac{1}{2\pi\sigma_z(t)\sigma_y(t)(1-\rho^2(t))^{1/2}}
\exp \left\{-\frac{1}{2(1-\rho^2(t))} \left[ \frac{(y-q^\epsilon(t))^2}{\sigma^2_y(t)}+\frac{(z-m^\epsilon(t))^2}{\sigma^2_z(t)}\right.\right.\nonumber\\
 & -  \left.\left.\frac{2\rho(t)(y-q^\epsilon(t))(z-m^\epsilon(t))}{\sigma_y(t)\sigma_z(t)}
\right] \right\},\label{chap4:pedens}
\end{align}
where the correlation coefficient $\rho(t)$ is defined by $\sigma_{yz}(t)/\sigma_{y}(t)\sigma_{z}(t)$.
Observe that for $\epsilon=0$, \eqref{chap4:pedens} reduces to
\begin{align}
p^0(y,z,t)=&\frac{1}{2\pi\sigma_z(t)\sigma_y(t)(1-\rho^2(t))^{1/2}}\exp \left\{-\frac{1}{2(1-\rho^2(t))} \left[ \frac{(y-q^0(t))^2}{\sigma^2_y(t)}+\frac{(z-m^0(t))^2}{\sigma^2_z(t)}\right.\right.\nonumber\\
 & - \left.\left.\frac{2\rho(t)(y-q^0(t))(z-m^0(t))}{\sigma_y(t)\sigma_z(t)} \right] \right\},
\label{chap4:ped0ns}
\end{align}
with
\begin{equation}
\label{chap4:meanvarz0}
m^0(t):=Ye^{-\frac{c_0t}{2}} (\cos \omega t+\frac{c_0}{2\omega} \sin \omega t ),\ \
q^0(t):=-Y\frac{k}{\omega}e^{-\frac{c_0t}{2}} \sin \omega t.
\end{equation}
From \eqref{chap4:meanvarze}-\eqref{chap4:meanvarz0} we can easily see that
\begin{equation}\label{chap4:metom0}
m^\epsilon(t)=m^0(t)-\epsilon f(t), \ \ \
q^\epsilon(t)=q^0(t)+\epsilon g(t),
\end{equation}
with
\[
f(t):=e^{-\frac{c_0t}{2}} (\cos \omega t+\frac{c_0}{2\omega} \sin \omega t ) \quad \mbox{and} \quad g(t):=\frac{k}{\omega} e^{-\frac{c_0t}{2}} \sin \omega t .
\]
Plugging \eqref{chap4:metom0} into \eqref{chap4:pedens}, we obtain
\begin{align}
p^\epsilon(y,z,t)=&\frac{1}{2\pi\sigma_z(t)\sigma_y(t)(1-\rho^2(t))^{1/2}}\exp \left\{-\frac{1}{2(1-\rho^2(t))} \left[ \frac{(y-\left[q^0(t)+\epsilon g(t)\right])^2}{\sigma^2_y(t)} +\right.\right.\nonumber\\
&   \left.\left. \frac{(z-\left[m^0(t)-\epsilon
f(t)\right])^2}{\sigma^2_z(t)} - \frac{2\rho(t)(y-\left[q^0(t)+\epsilon g(t)\right])(z-\left[m^0(t)-\epsilon f(t)\right])}{\sigma_y(t)\sigma_z(t)} \right] \right\}.\label{chap4:pedens2}
\end{align}
Now, notice that
\begin{align}
&\frac{(y-\left[q^0(t)+\epsilon g(t)\right])^2}{\sigma^2_y(t)}+\frac{(z-\left[m^0(t)-\epsilon f(t)\right])^2}{\sigma^2_z(t)}- \frac{2\rho(t)(y-\left[q^0(t)+\epsilon g(t)\right])(z-\left[m^0(t)-\epsilon f(t)\right])}{\sigma_y(t)\sigma_z(t)}\nonumber\\
&=\frac{(y-q^0(t))^2}{\sigma^2_y(t)}+\frac{(z-m^0(t))^2}{\sigma^2_z(t)} -\frac{2\rho(t)(y-q^0(t))(z-m^0(t))}{\sigma_y(t)\sigma_z(t)} +\epsilon^2\left[\frac{g^2(t)}{\sigma^2_y(t)}+\frac{f^2(t)}{\sigma^2_z(t)}+ \frac{2\rho(t)g(t)f(t)}{\sigma_y(t)\sigma_z(t)}\right]\nonumber\\
&-2\epsilon\left[\frac{(y-q^0(t))g(t)}{\sigma^2_y(t)}-\frac{(z-m^0(t))f(t)}{\sigma^2_z(t)}+ \frac{\rho(t)}{\sigma_y(t)\sigma_z(t)}\left((y-q^0(t))f(t)-(z-m^0(t))g(t)\right)\right].\label{chap4:expterm1}
\end{align}
Then considering \eqref{chap4:ped0ns}, we have
\begin{align}
p^\epsilon(y,z,t)=&p^0(y,z,t)\exp \left\{ -\frac{\epsilon^2}{2} A(t) + \epsilon \frac{ (y-q_0(t))r(t)-(z-m_0(t)) s(t)}{(1-\rho^2(t))\sigma_y(t)\sigma_z(t)}  \right \}.\nonumber\\
\label{chap4:pedens3}
\end{align}
where 
\begin{eqnarray*}
A(t)&:=&\frac{1}{1-\rho^2(t)}\left(\frac{g^2(t)}{\sigma^2_y(t)}+\frac{f^2(t)}{\sigma^2_z(t)}+
\frac{2\rho(t)g(t)f(t)}{\sigma_y(t)\sigma_z(t)}\right),\\
r(t) &:=& \frac{g(t)\sigma_z(t)}{\sigma_y(t)}+\rho(t)f(t),\\
s(t) &:=& \frac{f(t)\sigma_y(t)}{\sigma_z(t)}+\rho(t)g(t).\\
\end{eqnarray*}

We now give a representation of $u(0,Y-\epsilon,0)$ in terms of the densities $p^\epsilon $ and $p^0 $ of the Gaussian processes $(z_{0,Y-\epsilon}(t),y_{0,Y-\epsilon}(t))$ and  $(z_{0,Y}(t),y_{0,Y}(t))$ respectively.
The proof is postponed to Section \ref{chap4:prooflemma}.
\begin{lemma}\label{chap4:lemma1}
Let $u$ be a solution of \eqref{chap4:pde1}. Then, it satisfies
\begin{align}
u(0,Y-\epsilon,0)=&
\int_D\left[p^\epsilon(y,z,T)-p^0(y,z,T)\right] \d y \d z \nonumber\\
 + & \int_{D_T^-} yu(y,Y,t)p^0(y,Y,t)\left[\exp\left\{-\frac{1}{2}\epsilon^2A(t) +\frac{\epsilon (y r(t)-Yh(t))}{(1-\rho^2(t))\sigma_y(t)\sigma_z(t)}\right\}-1\right] \d y \d t\nonumber\\
  - & \int_{D_T^+} yu(y,-Y,t)p^0(y,-Y,t)\left[\exp\left\{-\frac{1}{2}\epsilon^2A(t) +\frac{\epsilon (y r(t)-Yl(t))}{(1-\rho^2(t))\sigma_y(t)\sigma_z(t)}\right\}-1\right]\d y \d t,\nonumber\\
\label{chap4:testpde4}
\end{align}
with
\begin{eqnarray*}
 Yh(t)&:=& q_0(t)r(t)+(Y-m_0(t)) s(t),  h(t):=-g(t)r(t)+(1-f(t))s(t),\\[2mm]
 Yl(t)&:=&q_0(t)r(t)-(Y+m_0(t))s(t),  l(t):=-g(t)r(t)-(1+f(t))s(t),\\[2mm]
 D_T^+&:=&(0,T) \times (0, \infty)\ {\rm and}\  D_T^-:=(0,T) \times (-\infty,0).
\end{eqnarray*}
\end{lemma}
Now consider the terms
\begin{align}\label{chap4:Ie}
H^\epsilon = & \int_D\left[p^\epsilon(y,z,T)-p^0(y,z,T)\right] \d y \d z, \nonumber\\
I^\epsilon = & \int_{D_T^-} yu(y,Y,t)p^0(y,Y,t)\left[\exp\left\{-\frac{1}{2}\epsilon^2A(t) +\frac{\epsilon (y r(t)-Yh(t))}{(1-\rho^2(t))\sigma_y(t)\sigma_z(t)}\right\}-1\right] \d y \d t,\\
J^\epsilon = & - \int_{D_T^+} yu(y,-Y,t)p^0(y,-Y,t)\left[\exp\left\{-\frac{1}{2}\epsilon^2A(t) +\frac{\epsilon (y r(t)-Yl(t))}{(1-\rho^2(t))\sigma_y(t)\sigma_z(t)}\right\}-1\right] \d y \d t. \nonumber
\end{align}
Next, we study the behavior of these last integrals, with $u$ satisfying \eqref{chap4:pde1} so that the previous lemma holds, when $\epsilon$ is sufficiently small. 
The proof is also postponed to Section \ref{chap4:prooflemma}.
\begin{lemma}\label{chap4:lemma2}
Let $J^\epsilon$, $I^\epsilon$, and $H^\epsilon$ be the integrals of above.  Suppose that 
\[
k > X_+(c_0) :=\frac{1}{2} \left ( -\frac{c_0}{3} + c_0\sqrt{\frac{1}{9} + 4 \frac{c_0}{6}} \right ).
\]
Then,
\begin{itemize}
\item
$\liminf_{\epsilon \to 0} \frac{I^\epsilon}{\epsilon} = +\infty$,
\item
$\lim_{\epsilon \to 0} \frac{J^\epsilon}{\epsilon}$ is finite,
\item
$\lim_{\epsilon \to 0} \frac{H^\epsilon}{\epsilon}$ is finite.
\end{itemize}
Therefore,
\[
\lim_{\epsilon \to 0} \frac{u(0,Y-\epsilon,0)}{\epsilon} = +\infty \quad \mbox{and} \quad \lim_{\epsilon \to 0} \frac{u(0,-Y+\epsilon,0)}{\epsilon} = +\infty.
\]
\end{lemma}

\subsection{Proof of Theorem \ref{chap4:thm}}
We shall use the notation $\sigma_n^\epsilon = \textup{sign}(z^\epsilon(\tau^\epsilon_{n}-))$. 
Recall that, for each $n\geq 1$, the stopping time $\tau_n^\epsilon$ represents the instant of the $n-$th jump of the process $(y^\epsilon(t),z^\epsilon(t))$. Hence, for all $\tau_{n}^\epsilon\leq t < \tau_{n+1}^\epsilon$ and $n\geq 1$, we deduce from \eqref{chap4:SVI} and \eqref{chap4:SVIe} that
\begin{eqnarray*}
&& \dot{y}(t)-\dot{y}^\epsilon(t)=-\left[c_0 (y(t)-y^\epsilon(t)) + k (z(t)-z^\epsilon(t))\right], \quad  \mbox{and}\\
&& (\dot{z}^\epsilon(t)-y^\epsilon(t))(z(t)-z^\epsilon(t))\geq 0,\\
&& (\dot{z}(t)-y(t))(z^\epsilon(t)-z(t))\geq 0.
\end{eqnarray*}
By using the notation $d/dt$ of derivatives, we obtain
\begin{eqnarray}
&& \frac{d}{dt}\left(y(t)-y^\epsilon(t)\right)=-c_0(y(t)-y^\epsilon(t)) - k(z(t)-z^\epsilon(t)),\label{chap4:rest-eq-l2}\\
&& \left(\frac{d}{dt}(z(t)-z^\epsilon(t))-(y(t)-y^\epsilon(t))\right)(z(t)-z^\epsilon(t))\leq 0,\label{chap4:rest-ineq-l2} 
\end{eqnarray}
Multiplying by $(y(t) - y^\epsilon(t))$ in \eqref{chap4:rest-eq-l2} and using the product rule for derivatives, we get from \eqref{chap4:rest-eq-l2} and \eqref{chap4:rest-ineq-l2}
\begin{eqnarray}
\frac{1}{2}\frac{d}{dt}\left \vert y(t)-y^\epsilon(t)\right \vert^2 + c_0 \vert y(t)-y^\epsilon(t) \vert^2 & \leq & -k(z(t)-z^\epsilon(t))(y(t)-y^\epsilon(t))\nonumber\\
 & \leq & -\frac{k}{2}\frac{d}{dt} \left \vert (z(t)-z^\epsilon(t))\right \vert^2 \label{chap4:bound-dif-l2-pre}
\end{eqnarray}
for all\ $\tau_{n}^\epsilon\leq t<\tau_{n+1}^\epsilon$ and $n\geq1$.
Now, integrating \eqref{chap4:bound-dif-l2-pre} on $[\tau_{n}^\epsilon,\tau_{n+1}^\epsilon)$ and noting that $y(\tau^\epsilon_n-)=y(\tau^\epsilon_n)$, $y^\epsilon(\tau^\epsilon_n-)=y^\epsilon(\tau^\epsilon_n)=0$, $z(\tau^\epsilon_n-)=z(\tau^\epsilon_n)$, for all $n\geq 1$, we obtain
\begin{equation}\label{chap4:e2.5}
\left \vert y(\tau^\epsilon_{n+1})\right \vert^2
-\left \vert y(\tau_{n}^\epsilon)\right \vert^2
+2c_0\int_{\tau^\epsilon_n}^{\tau^\epsilon_{n+1}}\left \vert y(t)-y^\epsilon(t)\right \vert^2\
\d t
+k\left \vert z(\tau^\epsilon_{n+1})-z^\epsilon(\tau^\epsilon_{n+1}-)\right \vert^2
-k\left \vert z(\tau_{n}^\epsilon)-z^\epsilon(\tau_{n}^\epsilon)\right \vert^2\leq
0.
\end{equation}
But
\begin{eqnarray}
k\left \vert z(\tau^\epsilon_{n})-z^\epsilon(\tau^\epsilon_{n})\right \vert^2
&=& k\left \vert \left(z(\tau^\epsilon_{n})-z^\epsilon(\tau^\epsilon_{n}-)\right) +\left(z^\epsilon(\tau^\epsilon_{n}-)-z^\epsilon(\tau^\epsilon_{n})\right)\right \vert^2\nonumber\\
&=& k\left \vert z(\tau^\epsilon_{n})-z^\epsilon(\tau^\epsilon_{n}-)\right \vert^2 + k\ \epsilon^2 \nonumber\\
&& +2k\epsilon\ \sigma_n^\epsilon \left(z(\tau^\epsilon_n)-z^\epsilon(\tau^\epsilon_{n}-)\right).\label{chap4:e2.6}
\end{eqnarray}
Plugging \eqref{chap4:e2.6} into \eqref{chap4:e2.5}, and rearranging terms, we obtain
\begin{eqnarray*}
&&\left \vert y(\tau^\epsilon_{n+1})\right \vert^2 -\left \vert y(\tau_{n}^\epsilon)\right \vert^2 +k\left \vert z(\tau^\epsilon_{n+1})-z^\epsilon(\tau^\epsilon_{n+1}-)\right \vert^2 -k\left \vert z(\tau_{n}^\epsilon)-z^\epsilon(\tau_{n}^\epsilon-)\right \vert^2 \nonumber\\
&&+2c_0\int_{\tau^\epsilon_n}^{\tau^\epsilon_{n+1}}\left \vert y(t)-y^\epsilon(t)\right \vert^2\ \d t \leq k \epsilon^2 + 2 k\epsilon (\sigma_n^\epsilon z(\tau_n^\epsilon)-Y ).
\label{chap4:bound-int-n+1-l2-fin-rewr}
\end{eqnarray*}
We can drop the term $2k \epsilon (\sigma_n^\epsilon z(\tau_n^\epsilon)-Y ) \leq 0$ and get
\begin{eqnarray}
&&\left \vert y(\tau^\epsilon_{n+1})\right \vert^2 -\left \vert y(\tau_{n}^\epsilon)\right \vert^2 +k\left \vert z(\tau^\epsilon_{n+1})-z^\epsilon(\tau^\epsilon_{n+1}-)\right \vert^2 -k\left \vert z(\tau_{n}^\epsilon)-z^\epsilon(\tau_{n}^\epsilon-)\right \vert^2\nonumber\\
&&+2c_0\int_{\tau^\epsilon_{n}}^{\tau^\epsilon_{n+1}}\left \vert y(t)-y^\epsilon(t)\right \vert^2\ \d t\leq k\epsilon^2.\label{chap4:e2.7}
\end{eqnarray}
Observe that, for $N \in \mathbb{N}^\star$, we can iterate  \eqref{chap4:e2.7} for $1 \leq n \leq N$ to obtain
\begin{eqnarray}
&&\left \vert y(\tau^\epsilon_{N+1})\right \vert^2 -\left \vert y(\tau_1^\epsilon)\right \vert^2 +k\left \vert z(\tau^\epsilon_{N+1})-z^\epsilon(\tau^\epsilon_{N+1}-)\right \vert^2 -k\left \vert z(\tau_1^\epsilon)-z^\epsilon(\tau_1^\epsilon-)\right \vert^2\nonumber\\
&&+2c_0\int_{\tau^\epsilon_1}^{\tau^\epsilon_{N+1}}\left \vert y(t)-y^\epsilon(t)\right \vert^2\ \d t\leq kN\epsilon^2. \nonumber
\end{eqnarray}
Also, recalling that $y(\tau^\epsilon_{1})=0$, $\left \vert z(\tau_{1}^\epsilon)-z^\epsilon(\tau_{1}^\epsilon-)\right \vert^2=0$, and that
$\int_{0}^{\tau^\epsilon_{1}}\left \vert y(t)-y^\epsilon(t)\right \vert^2 \d t =0$, we derive:
\begin{eqnarray}
&&\left \vert y(\tau^\epsilon_{N+1})\right \vert^2 +k\left \vert z(\tau^\epsilon_{N+1})-z^\epsilon(\tau^\epsilon_{N+1}-)\right \vert^2 +2c_0\int_{0}^{\tau^\epsilon_{N+1}}\left \vert y(t)-y^\epsilon(t)\right \vert^2\ \d t\leq k\epsilon^2N.\label{chap4:e2.8}
\end{eqnarray}
Denote the total number of jumps  of the process $(y^\epsilon(t),z^\epsilon(t))$ arising in the time interval $(0,T)$ by $N_T^\epsilon:=\max_N\left\{\tau_N^\epsilon\leq T\right\}$. Note that $T < \tau_{N_T^\epsilon +1}$. Hence, from \eqref{chap4:e2.8}, we deduce
\begin{equation}
\sup_{1\leq n\leq N_T^\epsilon+1}\left \vert y(\tau^\epsilon_{n})\right \vert^2 
+k\sup_{1\leq n\leq N_T^\epsilon+1}\left \vert z(\tau^\epsilon_{n})-z^\epsilon(\tau^\epsilon_{n}-)\right \vert^2 +2c_0\int_{0}^{T}\left \vert y(t)-y^\epsilon(t)\right \vert^2\ \d t\leq
k\epsilon^2N_T^\epsilon.\label{chap4:e2.9}
\end{equation}
Assume first $z(0)=Y-\epsilon$. According to the definition of \eqref{chap4:theta-eps} set $\theta^\epsilon:=\theta(0,Y-\epsilon) = \inf \{t>0, \quad \vert z^\epsilon(t) \vert =Y\}=\inf \{t>0, \quad \vert z_{0,Y-\varepsilon}(t) \vert =Y\}$. It is clear that $\tau_1^\epsilon>\theta^\epsilon$ a.s. and then $\mathbb{P}(\tau_1^\epsilon>T)>\mathbb{P}(\theta^\epsilon>T)$. Now, let us assume $z(0)=-Y+\epsilon$. It is easy to verify that $u(-y,-z,t)=u(y,z,t)$, which gives
\[
\mathbb{P}(\theta^\epsilon>T)=u(0,Y-\epsilon,0)=u(0,-Y+\epsilon,0).
\]
Thus, by Lemma \ref{chap4:lemma2} we have $\frac{\mathbb{P}(\theta^\epsilon>T)}{\epsilon}  \to +\infty$. Therefore, if the initial condition $z(0)$ associated to \eqref{chap4:SVIe}, is a random variable $\Gamma\overset{({\rm law})}{=} p_1\delta_{Y-\epsilon}+(1-p_1)\delta_{-Y+\epsilon} $ independent of the Wiener process $w(t)$, then again setting $\theta_\Gamma^\epsilon:=\inf\{t>0, \quad |z_{0,\Gamma}(t)|=Y \} $, 
\[
\frac{\mathbb{P}(\theta_\Gamma^\epsilon>T)}{\epsilon} \to +\infty \quad \mbox{as} \quad \epsilon \to 0.
\]
Coming back to equation \eqref{chap4:e2.8} and noting that $N_T^\epsilon = \sum_n \chi_{\{ \tau_n^\epsilon \leq T \}}$, we get
\begin{equation}
\label{chap4:NTindic}
\mathbb{E}N_T^\epsilon = \sum_{n=1}^\infty\mathbb{E}\chi_{\{\tau_n^\epsilon\leq T\}}=\mathbb{E}\chi_{\{\tau_1^\epsilon\leq T\}}+\sum_{n=2}^\infty \mathbb{E}\chi_{\{\tau_n^\epsilon\leq T\}}.
\end{equation}
Observe that for all $n\geq 2$ and that $\tau_{n}^\epsilon-\tau_{n-1}^\epsilon$ is independent of $\tau_{n-1}^\epsilon$.
\begin{equation}\label{chap4:NTindic2}
\mathbb{E}\chi_{\{\tau_n^\epsilon\leq T\}}=\mathbb{E}\left[\chi_{\{\tau_{n-1}^\epsilon\leq T\}} \chi_{\{\tau_{n}^\epsilon-\tau_{n-1}^\epsilon \leq T-\tau_{n-1}^\epsilon\}}\right]
\leq \mathbb{E}\chi_{\{\tau_{n-1}^\epsilon\leq T\}} \mathbb{E}\chi_{\{\tau_{n}^\epsilon-\tau_{n-1}^\epsilon\leq T\}}.
\end{equation}
But note that
\[
\mathbb{E}\chi_{\{\tau_{n}^\epsilon-\tau_{n-1}^\epsilon\leq T\}}\leq
\mathbb{P}(\theta^\epsilon\leq T).
\]
From the last inequality and using \eqref{chap4:NTindic2}, we deduce
\[\mathbb{E}\chi_{\{\tau_{n}^\epsilon\leq
T\}}\leq \mathbb{E}\chi_{\{\tau_{1}^\epsilon\leq
T\}}(1-u(0,Y-\epsilon,0))^{n-1}.
\]
This yields
\begin{equation}\label{chap4:NTindic3}
\mathbb{E}N_T^\epsilon\leq \mathbb{E}\chi_{\{\tau_{1}^\epsilon\leq
T\}}\frac{(1-u(0,Y-\epsilon,0))}{u(0,Y-\epsilon,0)}\leq
\frac{\epsilon \mathbb{E}\chi_{\{\tau_{1}^\epsilon\leq
T\}}}{\epsilon u(0,Y-\epsilon,0)}.
\end{equation}
Hence, from Lemmas \ref{chap4:lemma1} and \ref{chap4:lemma2}
\begin{equation}\label{chap4:NTindic4}
\epsilon\ \mathbb{E}N_T^\epsilon\to 0\ \ \ \mbox{as
$\epsilon\textcolor{black}{\to} 0$.}
\end{equation}
Thus, as $\epsilon$ goes to $0$, \eqref{chap4:e2.9} and \eqref{chap4:NTindic4} yield
\begin{equation}\label{chap4:bound-taun-ce}
\frac{1}{\epsilon}\left\{\mathbb{E}\left[\sup_{1\leq n\leq N_T^\epsilon+1}\left \vert y(\tau_n^\epsilon)\right \vert^2\right]
+2c_0\mathbb{E}\int_0^T\left \vert y(t)-y^\epsilon(t)\right \vert^2\ \d t
+k\mathbb{E}\left[\sup_{1\leq n\leq N_T^\epsilon+1}\left \vert z(\tau^\epsilon_n)-z^\epsilon(\tau^\epsilon_n-)\right \vert^2 \right]\right\}\to
0.\ \ \ 
\end{equation}
Since the forced jumps have magnitude $\epsilon$, this implies: 
\begin{equation}\label{chap4:bound-taun-ce2}
\frac{1}{\epsilon}\left\{\mathbb{E}\left[\sup_{1\leq n\leq
N_T^\epsilon+1}\left \vert y(\tau_n^\epsilon)\right \vert^2\right]
+2c_0\mathbb{E}\int_0^T\left \vert y(t)-y^\epsilon(t)\right \vert^2\ \d t
+k\mathbb{E}\left[\sup_{1\leq n\leq N_T^\epsilon+1}\left \vert z(\tau^\epsilon_n)-z^\epsilon(\tau^\epsilon_n)\right \vert^2 \right]\right\}\to
0.
\end{equation}
Also, by \eqref{chap4:bound-dif-l2-pre}, we can see that any
$\tau^\epsilon_n\leq t<\tau^\epsilon_{n+1}$ satisfies
\[
\vert y(t)-y^\epsilon(t) \vert^2-\vert y(\tau_n^\epsilon)\vert^2+k \vert z(t)-z^\epsilon(t) \vert^2-
k\left \vert z(\tau^\epsilon_n)-z^\epsilon(\tau^\epsilon_n)\right \vert^2\leq
0.
\]
This gives
\[
\sup_{\tau^\epsilon_n\leq t<\tau^\epsilon_{n+1}}
\left\{\vert y(t)-y^\epsilon(t)\vert^2+k \vert z(t)-z^\epsilon(t) \vert^2\right\} \leq
\vert y(\tau_n^\epsilon) \vert^2+k\left \vert z(\tau^\epsilon_n)-z^\epsilon(\tau^\epsilon_n)\right \vert^2.
\]
Hence,
\begin{align*}
\sup_{\tau^\epsilon_1\leq t<T}
\left\{\vert y(t)-y^\epsilon(t) \vert^2+k \vert z(t)-z^\epsilon(t) \vert^2\right\} &\leq&
\sup_{1\leq n\leq N_T^\epsilon+1}\left\{\sup_{\tau^\epsilon_n\leq
t<\tau^\epsilon_{n+1}}\left\{ \vert y(t) - y^\epsilon(t)\vert^2+k\left \vert z(t)-z^\epsilon(t)\right \vert^2\right\}\right\}\\
&&\leq \sup_{1\leq n\leq
N_T^\epsilon+1}\left\{ \vert y(\tau_n^\epsilon) \vert^2+k\left \vert z(\tau^\epsilon_n)-z^\epsilon(\tau^\epsilon_n)\right \vert^2\right\}.
\end{align*}
Also,
\[
\sup_{0\leq t\leq
T}\left\{\vert y(t)-y^\epsilon(t) \vert^2+k\left \vert z(t)-z^\epsilon(t)\right \vert^2\right\}\leq
\sup_{1\leq n\leq
N_T^\epsilon+1}\left\{\vert y(\tau_n^\epsilon) \vert^2+k\left \vert z(\tau^\epsilon_n)-z^\epsilon(\tau^\epsilon_n)\right \vert^2\right\}.
\]
Therefore, \eqref{chap4:bound-taun-ce2} gives
\[
\frac{1}{\epsilon}\mathbb{E}\left[\sup_{0\leq t\leq
T}\left\{\vert y(t)-y^\epsilon(t) \vert^2+k\left \vert z(t)-z^\epsilon(t)\right \vert^2\right\}\right]\to
0\ \ \ \mbox{as $\epsilon\textcolor{black}{\to} 0$.}
\]
\section{Proof of the technical lemmas}
\label{chap4:prooflemma}
This section is devoted to the proofs of Lemmas \ref{chap4:lemma1} and \ref{chap4:lemma2}.

{\itshape Proof of Lemma \ref{chap4:lemma1}.} From \eqref{chap4:pde1}, we have
\begin{eqnarray}
0&=&\int_0^T\int_D (-u_t+Au)p^\epsilon\ \d y \d z \d t\nonumber\\
&=&\int_0^T\int_D (-u_t-\dfrac{1}{2} u_{yy} + (c_0y+kz)  u_{y} - y u_z)p^\epsilon\ \d y \d z \d t\nonumber\\
&=& -\int_D p^\epsilon(y,z,T)dydz+u(0,Y-\epsilon,0)+\int_0^T\int_D up_t^\epsilon \d y \d z \d t\nonumber\\
&&-\int_0^T\int_D\frac{1}{2}up^\epsilon_{yy}dydzdt-\int_0^T\int_Du((c_0y+kz)p^\epsilon)_y \d y \d z \d t\nonumber\\
&&-\int_{D_T^-} yu(y,Y,t)p^\epsilon(y,Y,t)dydt+\int_{D_T^+} yu(y,-Y,t)p^\epsilon(y,-Y,t) \d y \d t\nonumber\\
&&+\int_0^T\int_Dyup^\epsilon_z \d y \d z \d t.\label{chap4:testpde1}
\end{eqnarray}
By using \eqref{chap4:bound-cond-pe} and rearranging terms, \eqref{chap4:testpde1} becomes
\begin{equation}
\label{chap4:testpde2}
u(0,Y-\epsilon,0)=\int_D p^\epsilon(y,z,T)\d y \d z+\int_{D_T^-} yu(y,Y,t)p^\epsilon(y,Y,t)\d y \d t -\int_{D_T^+}  yu(y,-Y,t)p^\epsilon(y,-Y,t)\d y \d t.
\end{equation}
In addition, $p^0(y,z,0):= \delta_{0,Y}(y,z)$, and
\begin{equation}
\label{chap4:dens0-prop}
0=\int_D p^0(y,z,T)\d y \d z+\int_{D_T^-}
yu(y,Y,t)p^0(y,Y,t) \d y \d t -\int_{D_T^+}
yu(y,-Y,t)p^0(y,-Y,t) \d y \d t.
\end{equation}
Using \eqref{chap4:pedens3} and substracting \eqref{chap4:dens0-prop} to \eqref{chap4:testpde2}, we can deduce the result \eqref{chap4:testpde4}.
\\\\
{\itshape Proof of Lemma \ref{chap4:lemma2}.} First note that, on a neighborhood of $t=0$, we have the following expansions: 

\begin{itemize}
\item
$f(t) = e^{-\frac{c_0t}{2}} (\cos \omega t+\frac{c_0}{2\omega} \sin \omega t ) = 1 - k \frac{t^2}{2} + \frac{c_0}{12} (c_0^2 + 2 \omega^2) t^3 + o(t^3)$,
\item
$g(t)= \frac{k}{\omega} e^{-\frac{c_0t}{2}} \sin \omega t = kt(1-\frac{c_0}{2}t)+o(t^2)$.
\end{itemize}
From \eqref{chap4:meanvarze}-\eqref{chap4:covyz}, we also have
\begin{itemize}
\item
$\sigma^2_y(t) = t-c_0t^2 + o(t^2), \quad \sigma_y(t) = \sqrt{t} \left ( (1-\frac{c_0}{2}t) + o(t) \right )$,
\item
$\sigma^2_z(t) = \frac{t^3}{3}-c_0\frac{t^4}{4} + o(t^4), \quad \sigma_z(t) = \frac{t^{\frac{3}{2}}}{\sqrt{3}} \left ( (1-\frac{c_03}{8}t) + o(t) \right )$,
\item
$\frac{\sigma_z(t)}{\sigma_y(t)} = \frac{t}{\sqrt{3}} \left ( (1+\frac{c_0}{8}t) + o(t) \right )$,
\item
$\rho(t) = \frac{\sqrt{3}}{2}\left ( 1-\frac{c_0}{8}t + o(t) \right )$, recalling that $\rho(t) = \frac{\sigma_{yz}(t)}{\sigma_y(t)\sigma_z(t)}$.
\end{itemize}
Equation \eqref{chap4:meanvarz0} yields
\begin{itemize}
\item
$q^0(t) = -Ykt(1-\frac{c_0}{2}t) + o(t^2)$,
\item
$m^0(t) = Y(1 - k \frac{t^2}{2}) + o(t^2)$.
\end{itemize} 
Recalling that $r(t) = \frac{g(t)\sigma_z(t)}{\sigma_y(t)}+\rho(t)f(t)$ and $s(t):=\frac{f(t)\sigma_y(t)}{\sigma_z(t)}+\rho(t)g(t)$ and using the previous estimations, we can check that $-g(t)r(t) = -\frac{k\sqrt{3}}{2} t + \frac{5\sqrt{3}}{16} c_0 k t^2 + o(t^2)$ and $(1-f(t))s(t)=\frac{k\sqrt{3}}{2} t  + 
\frac{\sqrt{3}}{4} \left ( k^2 -\frac{3}{4} c_0 k- \frac{c_0^3}{6} \right ) t^2 + o(t^2)$. Therefore,
 $h(t) \sim \frac{\sqrt{3}}{4} P(c_0,k) t^2$ where
\[
P(c_0,k) := k^2 + \frac{c_0}{3}k - \frac{c_0^3}{6}.
\]
Denote $X_+(c_0):= \frac{1}{2} \left ( -\frac{c_0}{3} + c_0\sqrt{\frac{1}{9} + 4 \frac{c_0}{6}} \right )$. Since we have assumed that $k > X_+(c_0)$,
 it then follows that $h^\prime(0)=0$ and  $h^{\prime \prime}(0)= \frac{\sqrt{3} P(c_0,k) }{2}>0$. \textcolor{black}{We can thus consider a fixed interval $(0,\tilde{t})$ such that $h^{\prime \prime}(t)>0$ on $[0,\tilde{t}]$, hence $h(t) > 0$ on $(0,\tilde{t})$. Also, we have $r(t)=\frac{\sqrt{3}}{2}(1-\frac{c_0}{8}t + o(t))$. 
 Hence, there exists a positive constant $\bar{t}$ such that $r(t)>0$ on $[0,\bar{t}]$.
Let $t_0:=\min\{\tilde{t},\bar{t}\}$. This implies that $h(t)>0$ and $r(t)>0$ on $(0,t_0)$.} Recall that $h(0)=0$. Now write, from \eqref{chap4:Ie}
\[
I^\epsilon=I^\epsilon_1+I^\epsilon_2,
\]
with
\begin{eqnarray*}
I^\epsilon_1&:=&\int_0^{t_0\wedge T}\int_{-\infty}^0 yu(y,Y,t)p^0(y,Y,t)\left[\exp\left\{-\frac{1}{2}\epsilon^2A(t) +\frac{\epsilon (y\textcolor{black}{r(t)}-Yh(t))}{(1-\rho^2(t))\sigma_y(t)\sigma_z(t)}\right\}-1\right] \d y \d t,\\
I^\epsilon_2&:=&\int_{t_0\wedge T}^T\int_{-\infty}^0 yu(y,Y,t)p^0(y,Y,t)\left[\exp\left\{-\frac{1}{2}\epsilon^2A(t) +\frac{\epsilon (y\textcolor{black}{r(t)}-Yh(t))}{(1-\rho^2(t))\sigma_y(t)\sigma_z(t)}\right\}-1\right] \d y \d t.
\end{eqnarray*}
From the definition of $t_0$, $h(t)\geq 0$ \textcolor{black}{and $r(t)>0$} for $0< t< t_0\wedge T$. \textcolor{black}{Moreover,} $y<0$ in $I^\epsilon_1$, so we have
\[
-\frac{1}{2}\epsilon^2A(t)+\frac{\epsilon(y\textcolor{black}{r(t)}-Yh(t))}{(1-\rho^2(t))\sigma_y(t)\sigma_z(t)}\leq 0.
\]
Therefore, the integrand in $I^\epsilon_1$ is a positive function. Now, using the basic inequality $\exp\{-x\}-1\leq -x\exp\{-x\}$,
for $x\geq 0$, we can write

\begin{eqnarray*}
\frac{I^\epsilon_1}{\epsilon}
& \geq & -\int_0^{t_0\wedge T}\int_{-\infty}^0 \left[yu(y,Y,t)p^0(y,Y,t)\left[\frac{1}{2}\epsilon A(t) -\frac{(y\textcolor{black}{r(t)}-Yh(t))}{(1-\rho^2(t))\sigma_y(t)\sigma_z(t)}\right]\right.\nonumber\\
&          & \times\left.\exp\left\{-\frac{1}{2}\epsilon^2A(t)+\frac{\epsilon(y\textcolor{black}{r(t)}-Yh(t))}{(1-\rho^2(t)) \sigma_y(t)\sigma_z(t)}\right\}\right] \d y \d t.
\end{eqnarray*}
As $A(t) \geq 0$, we get
\begin{eqnarray}\label{chap4:Ie1}
\frac{I^\epsilon_1}{\epsilon}
& \geq & \int_0^{t_0\wedge T}\int_{-\infty}^0 \left[yu(y,Y,t)p^0(y,Y,t)\left[ \frac{(y\textcolor{black}{r(t)}-Yh(t))}{(1-\rho^2(t))\sigma_y(t)\sigma_z(t)}\right]\right.\nonumber\\
&          & \times\left.\exp\left\{-\frac{1}{2}\epsilon^2A(t)+\frac{\epsilon(y\textcolor{black}{r(t)}-Yh(t))}{(1-\rho^2(t)) \sigma_y(t)\sigma_z(t)}\right\}\right] \d y \d t.
\end{eqnarray}
As the integrand in the right hand side of \eqref{chap4:Ie1} is a positive function, Fatou's lemma yields the following inequality,
\begin{eqnarray}\label{chap4:Ie1res}
\liminf_{\epsilon \to 0} \frac{I^\epsilon_1}{\epsilon}
& \geq & \int_0^{t_0\wedge T}\int_{-\infty}^0  yu(y,Y,t)p^0(y,Y,t)\left[ \frac{(y\textcolor{black}{r(t)}-Yh(t))}{(1-\rho^2(t))\sigma_y(t)\sigma_z(t)}\right]  \d y \d t.\nonumber\\
\end{eqnarray}
Note that in \refs{chap4:Ie1res} the right hand side may be $+\infty$.
For $I_2^\epsilon$, since $t\geq t_0\wedge T$, there is no singularity at $t=0$. Therefore, taking the limit of $I_2^\epsilon/\epsilon$, we obtain
\begin{equation}
\label{chap4:lim3}
\liminf_{\epsilon \to 0} \frac{I_2^\epsilon}{\epsilon} = \int_{t_0 \wedge T}^T\int_{-\infty}^0\frac{yu(y,Y,t)p^0(y,Y,t)(y\textcolor{black}{r(t)}-Yh(t))}{(1-\rho^2(t))\sigma_y(t)\sigma_z(t)}\d y \d t
\end{equation}
which is finite. 
Note that 
\[
J = - \int_0^T \frac{|h(t)|}{(1-\rho^2(t)) \sigma_y(t) \sigma_z(t)} \left [ \int_{-\infty}^0 y u(y,Y,t) p^0(y,Y,t) \d y \right ] \d t 
\]
is finite. Indeed, from the expansion of $h(t)$ we have that locally in time
\[
\frac{h(t)}{\sigma_y(t) \sigma_z(t)}
\]
is bounded. Moreover from \eqref{chap4:dens0-prop} above
\[
- \int_0^T \int_{-\infty}^0 y u(y,Y,t) p^0(y,Y,t) \d y  \d t < \infty.
\]
Collecting results we can assert that 
\begin{align}
\label{chap4:res}
\liminf_{\epsilon \to 0} \frac{I_\epsilon}{\epsilon} \geq & \int_0^T \int_{-\infty}^0 \frac{y^2\textcolor{black}{r(t)} u(y,Y,t) p^0(y,Y,t) \d y \d t}{(1-\rho^2(t)) \sigma_y(t) \sigma_z(t) }\\
&-\ Y\int_0^T \frac{h(t)}{(1-\rho^2(t))  \sigma_y(t)  \sigma_z(t)} \left( \int_{-\infty}^0 y u(y,Y,t) p^0(y,Y,t) \d y \right) \d t. \nonumber
\end{align}
The second integral is finite. Now, let us show that the first integral is $+\infty$.
We check that 
\[
\lim_{t \to 0} \int_{-\infty}^0 y^2 u(y,Y,t) p^0(y,Y,t) \d y > 0.
\]
The function $u(y,z,t)$ is increasing in $t$. Indeed, from the probabilistic representation we have 
\[
u(y,z,t_1) = \mathbb{P}[\theta(y,z)>T-t_1] \leq \mathbb{P}[\theta(y,z)>T-t_2] = u(y,z,t_2), \quad \forall t_1 \leq t_2.
\]
Therefore, we have
\begin{equation}
\label{chap4:lim4}
 \lim_{t \to 0} \int_{-\infty}^0 y^2 u(y,Y,0) p^0(y,Y,t) \d y  \leq \lim_{t \to 0} \int_{-\infty}^0 y^2 u(y,Y,t) p^0(y,Y,t) \d y
\end{equation}
Now, 
\[
u_y(0-,Y,0) <0. 
\]
Indeed, $\forall c >0, \quad u(-c,Y,0)>0$ and $u(0,Y,0)=0$. So,
\[
u_y(0-,Y,0)  \leq 0. 
\] 
It cannot be equal to $0$, otherwise the derivative exists and $u_y(0,Y,0)=0$. But then by minimum properties we have $u_{yy}(0,Y,0)>0$, that contradicts
\[
-u_t(0,Y,0) -\frac{1}{2}u_{yy}(0,Y,0) = 0.
\]
Therefore for $y<0$ close to $0$ we have 
\[
u(y,Y,0) \sim a y, \quad a <0.
\]
On the interval $(-\eta,0)$, we can assume $u(y,Y,0)>\frac{a}{2}y$.
So,
\[
\lim_{t \to 0} \int_{-\eta}^0 y^2 u(y,Y,0) p^0(y,Y,t) \d y \geq \lim_{t \to 0} \int_{-\eta}^0 \frac{a}{2} y^3 p^0(y,Y,t) \d y. 
\]
From \eqref{chap4:res}, since $\int_0^T \frac{r(t)}{(1-\rho^2(t))\sigma_y(t) \sigma_z(t)} \d t = +\infty$, it is sufficient to check the property
\[
\lim_{t \to 0} \int_{-\eta}^0 y^3 p^0(y,Y,t) \d y < 0.
\]
Set 
\[
\tilde{q}_0(t) := q_0(t) + \frac{\rho(t)Y(1-f(t)) \sigma_y(t)}{\sigma_z(t)},
\]
then
\[
p^0(y,Y,t) = \frac{1}{2 \pi \sigma_y(t) \sigma_z(t) (1-\rho^2(t))^{\frac{1}{2}}} \exp \left (-\frac{1}{2} \frac{Y^2(1-f(t))^2}{\sigma_z^2(t)} \right ) \exp \left (-\frac{(y-\tilde{q}_0(t))^2}{2(1-\rho^2(t))\sigma_y^2(t)} \right ).
\]
\textcolor{black}{Hence,} denoting 
\[
L_\eta := \int_{-\eta}^{0} \frac{y^3}{\sqrt{2\pi} \sigma_y(t)(1-\rho^2(t))^{1/2}} \exp \left (-\frac{(y-\tilde{q}_0(t))^2}{2(1-\rho^2(t)) \sigma_y^2(t)} \right ) \d y,
\]
we get
\begin{equation}\label{chap4:Leta1}
\int_{-\eta}^0 y^3 p^0(y,Y,t) \d y = \frac{1}{\sqrt{2\pi} \sigma_z(t)} \exp \left (-\frac{1}{2} \frac{Y^2 (1-f(t))^2}{\sigma_z^2(t)} \right ) L_\eta.
\end{equation}
In addition, by change of variables, we have 
\begin{equation}\label{chap4:Leta2}
L_\eta = \int_{\frac{-\eta-\tilde{q}_0(t)}{(1-\rho^2(t))^{1/2} \sigma_y(t)}}^{\frac{-\tilde{q}_0(t)}{(1-\rho^2(t))^{1/2} \sigma_y(t)}} (\tilde{q}_0(t) + (1-\rho^2(t))^{1/2} \sigma_y(t) u)^3 \exp(-\frac{1}{2} u^2) \frac{ \d u}{\sqrt{2 \pi}}.
\end{equation}
Therefore, for $t$ close to $0$ we have $\tilde{q}_0(t) \sim -\frac{Ykt}{4}$ and we can check using formula \eqref{chap4:Leta1} and \eqref{chap4:Leta2} that
\[
\lim_{t \to 0} \int_{-\eta}^0 y^3p^0(y,Y,t) \d y = \frac{1}{2 \pi} \frac{\sqrt{3}}{8} \int_{-\infty}^0 u^3 \exp(-\frac{1}{2} u^2) \d u.
\]
Finally, since $\frac{r(t)}{\sigma_y(t) \sigma_z(t)} = \frac{3}{2t^2} (1 +\frac{3 c_0}{4} t + o(t))$, the first integral in the right hand side in \eqref{chap4:res} is $+ \infty$.
We thus have proven
\[
\lim_{\epsilon \to 0} \frac{I^\epsilon}{\epsilon} = +\infty.
\]
Next, consider from \eqref{chap4:Ie} the term 
\[
\textcolor{black}{\frac{J^\epsilon}{\epsilon}=}-\frac{1}{\epsilon}\int_0^T\int_{0}^{\infty} y u(y,-Y,t)p^0(y,-Y,t)\left[\exp\left\{-\frac{1}{2}\epsilon^2A(t) +\frac{\epsilon (y\textcolor{black}{r(t)}-Yl(t))}{(1-\rho^2(t))\sigma_y(t)\sigma_z(t)}\right\}-1\right] \d y \d t.
\]
From \eqref{chap4:ped0ns} we have
\begin{align*}
p^0(y,-Y,t)
 = &\frac{1}{2\pi\sigma_z(t)\sigma_y(t)(1-\rho^2(t))^{1/2}}\exp \left\{-\frac{1}{2(1-\rho^2(t))} \left[ \frac{(y-q^0(t))^2}{\sigma^2_y(t)}+\frac{(Y+m^0(t))^2}{\sigma^2_z(t)}\right.\right.\nonumber\\
    & +\left.\left.\frac{2\rho(t)(y-q^0(t))(Y+m^0(t))}{\sigma_y(t)\sigma_z(t)} \right] \right\}.
\end{align*}
Due to the term $\exp\{-\frac{(Y+m^0(t))^2}{2(1-\rho^2(t))\sigma^2_z(t)}\}$, we do not have a singularity because of $\sigma_z(t)$. 
Indeed, $Y + m^0(t) \geq Y (1-\exp(-\frac{c_0 \pi}{2 \omega}))$, and for $t$ close to $0$ we have
\begin{align*}
& \frac{1}{\sigma_z(t)} \exp \left ( -\frac{(Y+m^0(t))^2}{2(1-\rho^2(t)) \sigma_z^2(t)} \right )\\  
& \leq   \frac{\sqrt{2(1-\rho^2(t))}}{Y (1-\exp(-\frac{c_0 \pi}{2 \omega}))} \left ( \frac{Y+m^0(t)}{\sqrt{2(1-\rho^2(t))} \sigma_z(t)} \right )  \exp \left ( -\frac{(Y+m^0(t))^2}{2(1-\rho^2(t)) \sigma_z^2(t)} \right ) \nonumber \\
& \leq  C \exp \left ( -\frac{(Y+m^0(t))^2}{4(1-\rho^2(t)) \sigma_z^2(t)} \right ) \nonumber
\end{align*}
where $C>0$ is a constant depending on $T$. From the above equation and recalling the asymptotics $\sigma_y(t) \sim \sqrt{t}$ for $t$ close to $0$, we deduce that the quantity
\[
-\int_0^T\int_{0}^{\infty} yu(y,-Y,t)p^0(y,-Y,t) [ \frac{(y\textcolor{black}{r(t)}-Yl(t))}{(1-\rho^2(t))\sigma_y(t)\sigma_z(t)}] dydt
\]
is well defined. In the same way
\[
\textcolor{black}{\frac{H^\epsilon}{\epsilon}=}\frac{1}{\epsilon}\int_D\left(p^\epsilon(y,z,T)-p^0(y,z,T)\right) \d z \d y,
\]
has a well defined limit. \textcolor{black}{Therefore, from \eqref{chap4:testpde4}, we deduce}
\begin{equation}
\label{chap4:limfin}
\frac{u(0,Y-\epsilon,0)}{\epsilon}\to +\infty \quad \mbox{as} \quad \epsilon\textcolor{black}{\to} 0.
\end{equation}
\textcolor{black}{As before, let us assume that $z(0)=-Y+\epsilon$. It is easy to see that $u(-y,-z,t)=u(y,z,t)$. This yields 
\[
\mathbb{P}(\theta^\epsilon>T)=u(0,Y-\epsilon,t)=u(0,-Y+\epsilon,0),
\]
so $\frac{u(0,-Y+\epsilon,0)}{\epsilon}\to +\infty$ as well. This completes the proof.}

\end{document}